\documentclass[11pt]{amsart}
\usepackage{amssymb}
\usepackage{palatino}
\usepackage{tikz-cd}
\input amssym.def

\usepackage{amsmath, amsfonts}
\usepackage{amssymb}
\usepackage{amscd}
\usepackage[mathscr]{eucal}
\usepackage{palatino}
\usepackage{comment}
\usepackage{enumerate}
\usepackage{enumitem}
\usepackage{xcolor}

\newfont{\cyrr}{wncyr10}

\newcommand{\thmref}[1]{Theorem~\ref{#1}}

\newcommand{\propref}[1]{Proposition~\ref{#1}}
\newcommand{\lemref}[1]{Lemma~\ref{#1}}

\newtheorem{thm}{Theorem}

\newtheorem{lem}[thm]{Lemma}

\newtheorem{prop}[thm]{Proposition}

\theoremstyle{remark}
\newtheorem{rmk}{Remark}
\newtheorem{quest}[rmk]{Question}

\def\({\left(}
\def\){\right)}
\def\[{\left[}
\def\]{\right]}

\def\N{\mathbb{N}}
\def\Z{\mathbb{Z}}

\def\cA{\mathcal{A}}
\def\cB{\mathcal{B}}
\def\cC{\mathcal{C}}

\def\cP{\mathcal{P}}

\def\cS{\mathcal{S}}

\newcommand{\lf}{\lfloor}
\newcommand{\rf}{\rfloor}
\renewcommand{\mod}{ \text{ mod }}

\title{Coprimality of elements in regular sequences with polynomial growth}

\author{Jean-Marc Deshouillers and Sunil Naik}

\address{Jean-Marc Deshouillers \newline
	Institut de Math\'ematiques de Bordeaux,
	Universit\'e de Bordeaux, CNRS, Bordeaux INP
	33400, Talence, France}
\email{jean-marc.deshouillers@math.u-bordeaux.fr}

\address{Sunil Naik \newline
	Department of Mathematics,
	Queen's University, Jeffrey Hall, 
	99 University Avenue, 
	Kingston, ON K7L3N6, 
	Canada}
\email{naik.s@queensu.ca}

\begin{document}
	
	\hfuzz 5pt
	
\subjclass[2020]{11B05, 11B25 11B50, 11K31, 11N56, 41A58}

\keywords{Regular sequences, $k$-coprimality, Banach density}
	
\maketitle
	
\begin{flushright}
{\em To Krishnaswami Alladi for his contribution to mathematics and to its diffusion}
\end{flushright}	
	
\begin{abstract}
The investigation of primes in certain arithmetic sequences is 
one of the fundamental problems in number theory and especially, 
finding blocks of distinct primes has gained a lot of attention in recent years. 
In this context, we prove the existence of long blocks of $k$-wise coprime elements in certain regular sequences. 
More precisely, we prove that for any positive integers $H \geq k \geq 2$ and 
for a real-valued $k$-times continuously differentiable  function 
$f \in \mathcal{C}^k\left( [1, \infty)\right)$ satisfying 
$\lim_{x \to \infty} f^{(k)}(x) = 0$ and $\limsup_{x \to \infty} f^{(k-1)}(x) = \infty$, 
there exist infinitely many positive integers $n$ such that
$$
\gcd\left( \lfloor f(n+i_1)\rfloor, \lfloor f(n+i_2)\rfloor, \cdots, \lfloor f(n+i_k)\rfloor \right) ~=~ 1
$$
for any integers $1 \leq i_1 < i_2  < \cdots < i_k \leq H$. Further, 
we show that there exists a subset $\cA \subseteq \N$ having upper Banach density one such that 
$$
\gcd\(\lf f(n_1) \rf, \lf f(n_2) \rf, \cdots, \lf f(n_k) \rf\) ~=~ 1
$$
for any distinct integers $n_1, n_2, \cdots, n_k \in \cA$. 
\end{abstract} 

\section{Introduction}
G. Lejeune Dirichlet \cite{Dir} proved in 1849  the following result
\begin{equation}\label{dir}
	\lim_{N \to \infty} \frac{1}{N^2}\#\{(n, m)\colon 1 \le n, m \le N, \gcd(n,m) =1\} ~=~ \frac{1}{\zeta(2)},
\end{equation}
where $\zeta$ denotes the Riemann zeta function. This question was popularized by Chebyshev who was often asking ``{}What is the probability that a given \textit{random} fraction $a/b$ with positive integers $a$ and $b$ is irreducible?". 
We outline the proof given by F. Mertens \cite{M} in 1874, which uses the M\"{o}bius function to detect that the $gcd$ of $n$ and $m$ is equal to $1$. We have
\begin{eqnarray*}
	\sum_{\substack{n, m \le N \\ \gcd(n,m)=1}}1 & = & \sum_{n, m \le N} \sum_{d \mid \gcd(n, m)} \mu(d)\\
	&=& \sum_{d \le N}\mu(d)\sum_{\substack{n \le N\\ n \equiv 0 \, (\text{mod } d)}}1 \sum_{\substack{m \le N\\ m \equiv 0 \,(\text{mod } d)}}1\\
	&=& \sum_{d \le N}\mu(d)\left(\frac{N}{d}+O(1)\right)^2 \sim \frac{N^2}{\zeta(2)}.
\end{eqnarray*}
Let $\cA=(a_n)_n$ and $\cB=(b_m)_m$ be two sequences of positive integers.  The question of the probability that two elements from $\cA$ and $\cB$ are coprime can be approached in two ways. We may first look for the asymptotic frequency with which $a_n$ and $b_m$ are coprime: Mertens' approach tells us that this amounts to knowing the distribution of the elements of $\cA$ and that of $\cB$ in sets of multiples. A second way is to look for the asymptotic frequency with which $a_n$ and $b_n$ are coprime: Mertens' approach tells us that this amounts to knowing the joint distribution of the pairs $(a_n, b_n)$ in sets of multiple, more precisely, we are reduced to studying expressions like 
\begin{equation}\label{mertens}
	\sum_{\substack{n \le N\\ a_n \equiv b_n \equiv 0 \, (\text{mod } d)}} 1
\end{equation}
We are only considering questions of the second type.\\

G. Watson \cite{Wa}, answering a question of K. F. Roth, proved in 1953 that for any irrational real number $\alpha$, the \textit{probability}\footnote{We say that the set of the integers satisfying a property $\cP$ has \textit{probability} $\delta$ if this set has \textit{natural asymptotic density} $\delta$, i.e. if $\lim_{N \to \infty}\#\{n \le N \colon \cP(n)\}/N = \delta$.} that $n$ and $\lf \alpha n \rf$ are coprime is $1/\zeta(2)$. J. Lambek and L. Moser \cite{LM} considered in 1955 the coprimality of $n$ and $\lf f(n) \rf$, when $f$ is a smooth slowly increasing function; this was extended and refined by P. Erd\H{o}s and G. G. Lorentz \cite{EL} in 1958. In 2002, F. Delmer and the first-named author \cite{DD} showed that for any non-integral $c$ larger than $1$, the \textit{probability} that $n$ and $\lf n^c \rf$ are coprime is still $1/\zeta(2)$, resolving a question formulated at the occasion of the workshop ``{}Paul Erd\H{o}s and his mathematics", held in Budapest in 1999. In 2017, V. Bergelson and F. K. Richter \cite{BR} extend the latter result by proving that for functions $f_1, \ldots, f_k$ belonging to a Hardy field and satisfying some growth conditions, the set of positive integers $n$ such that $\gcd(n, \lf f_1(n) \rf, \ldots,  \lf f_k(n) \rf)=1$ has \textit{probability} $1/\zeta(k+1)$. \\

In all the cases of the previous chapter, one of the sequence under consideration is the set of all the integers. This has two consequences on the sum (\ref{mertens}) : first, the numbers $d$ under consideration are bounded by $N$; second, the condition $a_n \equiv 0 \, (\text{mod } d)$ is easy to handle. Motivated by this consideration, M. Drmota, C. M\"{u}llner and the first-named author \cite{DDM} proved in 2023 the following theorem.
\begin{thm}\label{thmDDM}
	Let $c$ be in $(1, 2)$ and $\alpha < (c-1)$. There exist infinitely many $n$ such that, for any integer $H \le \alpha \log n$, all the elements in the sequence $\{\lfloor n^c \rfloor, \lfloor (n+1)^c \rfloor, \ldots, \lfloor (n+H)^c \rfloor\}$ are pairwise coprime.
\end{thm}
In \cite{DN}, we extended the qualitative part of Theorem \ref{thmDDM}, proving the following result.
\begin{thm}\label{thmDN1}
	Let $f$ be a twice continuously differentiable 
	real valued function defined on $[1, \infty)$ such that
	\begin{eqnarray}\label{2zeroatinf}
		& & \lim_{x \to \infty} f^{''}(x) ~=~ 0 \\
		\notag
		& & \text{and }\\
		\label{2unbounded}
		& & \limsup_{x \to \infty} f'(x) = \infty.
	\end{eqnarray}
	Then for any positive integer $H$, 
	there exist infinitely many positive integers $n$ such that 
	all the elements in the sequence 
	$\{ \lf f(n+1) \rf,  \lf f(n+2) \rf, \cdots, \lf f(n+H) \rf \}$ 
	are pairwise coprime.
\end{thm}
The main result of this paper is an extension of Theorem \ref{thmDN1} to more quickly increasing functions, namely functions with polynomial growth.

\begin{thm}\label{thmmain}
Let $k \geq 2$ be a positive integer and $f $  be a $k$-continuously differentiable real-valued function defined on $[1, \infty)$ such that
	\begin{eqnarray}\label{zeroatinf}
		& & \lim_{x \to \infty} f^{(k)}(x) ~=~ 0 \\
		\notag
		& & \text{and }\\
		\label{plusatinf}
		& & \limsup_{x \to \infty} f^{(k-1)}(x) = \infty.
	\end{eqnarray}
	Then for any positive integer $H\geq k$, there exist infinitely many 
	positive integers $n$ such that
	\begin{equation}\label{gcd1}
		\forall  \, h_1 < h_2 < \cdots < h_k \in [1, H] \colon \gcd\(\lf f(n + h_1) \rf, \lf f(n + h_2) \rf, \cdots, \lf f(n + h_k) \rf \) ~=~ 1.
	\end{equation}
\end{thm}

\begin{rmk}
	Condition (\ref{zeroatinf}) is not sufficient to insure the validity of the conclusion. For example, the function $f(x)=x^{k-1}+1/x$ satisfies (\ref{zeroatinf}) but for any positive $n$, two of the numbers $\lf f(n + 1) \rf, \lf f(n + 2) \rf, \lf f(n + 3) \rf, \lf f(n + 4) \rf$ are even.
\end{rmk}

\begin{rmk}
	With similar argument, we can prove Theorem \ref{thmmain} replacing Condition (\ref{plusatinf}) by the weaker requirement
	\begin{equation}\label{unbounded}
		f^{(k-1)} \text{ is unbounded}.
	\end{equation}
	One has to be slightly careful : if $x$ is a positive number which is not an integer, then $\lf -x \rf = -\lf x \rf -1$.
\end{rmk}

\begin{quest}
	Is it true that, with the hypotheses of Theorem \ref{thmmain}, there are integers $n$ such that the numbers $\lf f(n + 1) \rf, \lf f(n + 2) \rf, \ldots, \lf f(n + H) \rf$ are pairwise coprime?
\end{quest}

Our approach is based on using Taylor's expansion in the following way. Let $n$ and $u$ be two positive integers; there exists $\theta$ in $[0,1]$ such that  
\begin{eqnarray}
	\label{integralpart}
	f(n+u) &=& \sum_{j=0}^{k-1} \frac{u^j}{j!} \lf f^{(j)}(n)\rf\\
	\label{smallpart}
	&+& \{ f(n)\}\\
	\label{errorterm}
	&+& \sum_{j=1}^{k-1} \frac{u^j}{j!} \{ f^{(j)}(n)\} + \frac{u^k}{k!} f^{(k)}(n+\theta u).
\end{eqnarray}

The only knowledge we have on the term in (\ref{errorterm}) comes from (\ref{zeroatinf}): it will be treated as an error term having a small absolute value, but we'll have to keep in mind that we have no control on its sign. 

The main term is the RHS of (\ref{integralpart}): it is a rational number which is easily turned into an integer, on which we only need to know its residue in a class the modulus of which depends on $k$ and $H$ only. 
%This may come directly from the knowledge of the terms $\lf f^{(j)}(n)\rf$, with the help of $u$. But it may need to be slightly modified, in which case we shall need the term in (\ref{smallpart}).

In the case when the integral term has the expected arithmetic property, we expect to have 
\begin{equation}\label{bestcase}
	\lf f(n+u) \rf = \sum_{j=0}^{k-1} \frac{u^j}{j!} \lf f^{(j)}(n)\rf.
\end{equation}
For this, it is enough to have the \textit{positive} term in (\ref{smallpart}) small enough to be an error term but large enough to compensate the other error term in (\ref{errorterm}) in case it were negative. This will be illustrated in the proof of Proposition \ref{reduction}. In the case when the integral term needs some tuning, we shall also use the term in (\ref{smallpart}) as well as the first term in (\ref{errorterm}) to perform it. This will be illustrated in the proof of Theorem \ref{thmmain}.

The next proposition is an intermediate step in the proof of Theorem \ref{thmmain}: it gives sufficient conditions on the derivatives of $f$ at an integer $n$ for (\ref{gcd1}) to hold.
\begin{prop}\label{reduction}
	Let $H, k$ be positive integers with $2 \leq k \leq H$  and 
	$\Pi_H = \prod_{p \leq H} p$ be the product of all primes less than or equal to $H$. 
	Let $f \in \cC^k([1, \infty))$ be a real-valued function and $n$ be a positive integer such that
	\begin{equation}\label{c0pr}
		\frac{1}{4^{4k}} < \{f(n)\} < \frac{1}{4}, 
	\end{equation}
	\begin{equation}\label{c1pr}
		\{f^{(i)}(n)\} ~<~ \frac{1}{4H^i} ~\text{ for }~ 1 \leq i < k ~\text{ and }~
		|f^{(k)}(x)| ~<~ \frac{1}{{4^{4k}}H^k} ~\text{ for } x \in [n, n+H],
	\end{equation}
	\begin{equation}\label{c2pr}
		\forall i \in [1, k-1] \colon \lf f^{(i)}(n) \rf   ~\equiv~ 0 ~(\mod k! \Pi_H), \phantom{m}
	\end{equation}
	\begin{equation}\label{c3pr}
		\exists \, \ell \in [1, k-1] \colon \gcd\(\lf f(n) \rf,~ \lf f^{(\ell)}(n) \rf\) ~=~ 1.
	\end{equation}
	Then we have
	$$
	\gcd\( \lf f(n + h_1)\rf,  \lf f(n + h_2)\rf, \cdots,  \lf f(n + h_k)\rf\) ~=~ 1 
	$$
	for any integers $1 \leq h_1 < h_2 < \cdots < h_k \leq H$.
\end{prop}

We shall prove Proposition \ref{reduction} in Section \ref{secreduction} and apply it to prove Theorem \ref{thmmain} in Section \ref{Sthmmain}.\\

In Section 4, we shall extend Theorem \ref{thmmain} to an infinite set of $k$-wise coprime elements from $(\lf f(n) \rf)_n$, namely
\begin{thm}\label{thmB}
	Let $f$ be as in Theorem \ref{thmmain}. There exists an infinite subset $\cA \subset \mathbb{N}$ having upper Banach density equal to $1$ such that 
	$$
	\gcd\left(\left\lf f(n_1) \right\rf, \left\lf f(n_2) \right\rf, \ldots, \left\lf f(n_k) \right\rf\right) =1
	$$
	for any set of $k$ distinct integers in $\cA$.
\end{thm}
We recall that the upper Banach density (cf.  \cite{GTT} or \cite{Ri}) of a subset $\cS$ of $\mathbb{N}$ is defined by
$$
\lim_{H \to \infty} \limsup_{x \to \infty} \frac{1}{H}\#(\cS \cap (x, x+H]).
$$

\medspace

\section{PROOF OF PROPOSITION \ref{reduction}}\label{secreduction}

Let $f, n$ be as in Proposition \ref{reduction} and  $1 \leq h \leq H$. By Taylor's theorem,
we have for some $\theta_h$ in $[0,1]$:
\begin{eqnarray}
	\label{intP2}
	f(n+h) &=& \sum_{i=0}^{k-1} \frac{h^i}{i!} \lf f^{(i)}(n) \rf \\
	\label{fracP2}
	&+& \{f(n)\}\\
	\label{errP2}
	&+& \sum_{i=1}^{k-1} \frac{h^i}{i!} \{ f^{(i)}(n) \} 
	~+~ \frac{h^k}{k!} f^{(k)}(n + \theta_h h).
\end{eqnarray}

We have
$$
0 \le \frac{1}{4^{4k}} - \frac{h^k}{k! 4^{4k}H^k }\le \{f(n)\} + \sum_{i=1}^{k-1} \frac{h^i}{i!} \{ f^{(i)}(n) \} 
~+~ \frac{h^k}{k!} f^{(k)}(n + \theta_h h) \le \sum_{i=0}^{\infty} \frac{h^{i}}{i!}\cdot \frac{1}{4H^{i}} \le \frac{e}{4} < 1.
$$
By (\ref{c2pr}), for any $i$ between $0$ and $k-1$, the quantity $\lf f^{(i)}(n) \rf$ is divisible by $i!$. We thus have
\begin{equation}\label{eqIntf(n+h)}
	\forall h \in [1, H] \colon \lf f(n + h) \rf ~=~ \sum_{i = 0}^{k-1} \frac{h^i}{i!} \lf f^{(i)}(n) \rf .
\end{equation}

Let $ 1 \leq h_1 < h_2 < \cdots < h_k \leq H$ be positive integers.

Observe that for any prime $p \leq H$, we have
$$
p \nmid \lf f(n + h_j) \rf \phantom{m}\forall~ 1\le j \le k.
$$
This is because if $p \le H$, then from \eqref{c2pr}, 
we get $p \mid \frac{\lf f^{(i)}(n) \rf}{i!}$ for $1 \leq i \leq k-1$ and 
from \eqref{c3pr}, we have $p \nmid \lf f(n) \rf$.
Suppose that there exists a prime $p > H$ such that
\begin{equation}\label{eqp|gcd}
	p \mid \gcd\( \lf f(n + h_1)\rf,  \lf f(n + h_2)\rf, \cdots,  \lf f(n + h_k)\rf\).
\end{equation}
Then we rewrite \eqref{eqIntf(n+h)} in the matrix form as
\begin{equation}\label{eqVanMf}
	\begin{pmatrix}
		\lf f(n + h_1) \rf \\
		\lf f(n + h_2) \rf \\
		\vdots \\
		\lf f(n + h_k) \rf
	\end{pmatrix}
	~=~ 
	\begin{pmatrix}
		1 & h_1 & h_1^2 & \cdots & h_1^{k-1} \\
		1 & h_2 & h_2^2 & \cdots & h_2^{k-1} \\
		\vdots & \vdots & \vdots & \vdots & \vdots \\
		1 & h_k & h_k^2 & \cdots & h_k^{k-1} 
	\end{pmatrix}
	\begin{pmatrix}
		\lf f(n) \rf \\
		\frac{\lf f'(n) \rf }{1!}  \\
		\vdots\\
		\frac{\lf f^{(k-1)}(n) \rf }{(k-1)!}
	\end{pmatrix}
	.
\end{equation}
Let $V = V(h_1, h_2, \cdots, h_k) = (h_i^{j-1})_{1 \leq i, j \leq k}$ denote the Vandermonde matrix and let $V^{-1} = (b_{ij})_{1 \leq i, j \leq k}$ denote its inverse. It is well-known that
\begin{equation*}
	\det(V) ~=~ \prod_{1 \leq i < j \leq k} (h_j - h_i).
\end{equation*}
Also we have $\det(V) \cdot b_{ij} \in  \Z$ for every $1 \leq i, j \leq k$.
From \eqref{eqVanMf}, we have 
\begin{equation*}
	\frac{\lf f^{(r)}(n)\rf}{r!} ~=~
	\sum_{j = 1}^{k} b_{rj} \lf f(n + h_j) \rf
\end{equation*}
for any $0 \leq r \leq k-1$. Thus we get
\begin{equation*}
	\det(V) \cdot \lf f^{(r)}(n)\rf ~=~ r! \sum_{j = 1}^{k} \det(V) b_{rj} \lf f(n + h_j) \rf
\end{equation*}
From \eqref{eqp|gcd}, we get $p \mid \det(V) \cdot \lf f^{(r)}(n)\rf$ for any $0 \leq r \leq k-1$. Since $p \nmid \det(V)$, we get
$p \mid  \lf f^{(r)}(n)\rf$ for any $0 \leq r \leq k-1$, a contradiction to \eqref{c3pr}. This completes the proof of Proposition \ref{reduction}. \qed

\medspace

\section{Proof of \thmref{thmmain}} \label{Sthmmain}

The following result is a crucial step in proving \thmref{thmmain} by induction. It shows that if one can localize $\lf f^{(i)}(n)\rf$ in some congruence and $\{ f^{(i)}(n)\}$ in a small interval for some integer $n=n_0$ and $i$ in $[m, k-1]$, then we can do the same for $n=n_1$ close to $n_0$ and all $i$ in $[m-1, k-1]$.

\begin{lem}\label{thmInd}
	Let $k,m$ be positive integers with $m \leq k-1$ and 
	$\(A_i\)_{i=m}^{k-1}$, $\(B_i\)_{i = m}^{k}$ be strictly increasing sequence 
	of positive integers.  Also let $f \in \cC^{k}\([1, \infty)\)$ be a real-valued function 
	and $(v_i)_{i=m}^{k-1}$ and $n_0$ be positive integers such that
	\begin{equation}\label{cn0fInd}
		\lf f^{(i)}(n_0) \rf ~\equiv~ v_i ~(\mod k!\Pi_H) 
		\phantom{m}\text{ and }\phantom{m}
		\{f^{(i)}(n_0)\} ~\in~ \(\frac{1}{A_i}, \frac{1}{B_i}\),
		\phantom{m} 
		m \leq i \leq k-1
	\end{equation}
	and
	\begin{equation}\label{cfkInd}
		|f^{(k)}(x)| ~<~ \frac{1}{B_k} \phantom{m}\text{for } x \geq n_0.
	\end{equation}
	Suppose that 
	\begin{equation}\label{cBi6Ind}
		B_i ~<~ A_i \phantom{m}\text{ for } m \leq i \leq k-1, 
		\phantom{mm} 
		B_i ~>~ B_{i-1}^6 \phantom{m}\text{ for } m+1 \leq i \leq k,
	\end{equation}
	\begin{equation}\label{cB_kInd}
		A_m ~<~ \frac{B_m^2}{16k! \Pi_H}
		\phantom{m}\text{and}\phantom{m}
		B_k ~>~ \(4 k! \Pi_H A_{k-1}\)^{k+1}.
	\end{equation}
	Then for any integer $v_{m-1} \in \Z$ and positive integers  $B_{m-1} < A_{m-1} $ satisfying
	\begin{equation}\label{cAm-1Bm-1Ind}
		20 \frac{B_{m-1}}{A_{m-1}} ~<~ \frac{B_{m}}{A_{m}} 
		\phantom{m}\text{ and }\phantom{m}
		k! \Pi_H A_{m-1} ~<~ A_m,
	\end{equation}
	there exists a positive integer 
	$n_1 \in [n_0, n_0 + 4k!\Pi_H A_m]$ such that
	\begin{equation*} 
		\lf f^{(i)}(n_1) \rf ~\equiv~ v_i ~(\mod k!\Pi_H) 
		\phantom{m}\text{ and }\phantom{m}
		\{f^{(i)}(n_1)\} ~\in~ \(\frac{1}{2A_i}, \frac{2}{B_i}\)
	\end{equation*}
	for $m-1 \leq i \leq k-1$.
\end{lem}

\subsection*{Proof of \lemref{thmInd}}
Let the notations be as in \lemref{thmInd}. 
By Taylor's theorem, for any non-negative integer $h$, we have
\begin{equation*}
	f^{(m-1)}(n_0 + h) ~=~ \sum_{i = 0}^{k-m} \frac{h^i}{i!} f^{(m+i-1)}(n_0) 
	~+~ \frac{h^{k-m+1}}{(k-m+1)!} f^{(k)}(n_0+\theta_h h)
\end{equation*}
for some $\theta_h \in [0,1]$. 
Set
\begin{equation}\label{eqahdef}
	\begin{split}
		a_h 
		&~=~ f^{(m-1)}(n_0+h) 
		~-~  \sum_{i = 1}^{k-m} \frac{h^i}{i!} \lf f^{(m+i-1)}(n_0) \rf \\
		&~=~ f^{(m-1)}(n_0) ~+~ 
		\sum_{i = 1}^{k-m} \frac{h^i}{i!} \big\{ f^{(m+i-1)}(n_0) \big\} 
		~+~ \frac{h^{k-m+1}}{(k-m+1)!} f^{(k)}(n_0+\theta_h h).
	\end{split}
\end{equation}
Note that 
\begin{equation}\label{eqah+1-ah}
	\begin{split}
		a_{h+1} - a_h 
		~=~& \sum_{i=1}^{k-m} \frac{(h+1)^i - h^i}{i!} \big\{ f^{(m+i-1)}(n_0) \big\} 
		~+~ \frac{(h+1)^{k-m+1}}{(k-m+1)!} f^{(k)}(n_0+\theta_{(h+1)} (h+1)) \\
		&~-~ \frac{h^{k-m+1}}{(k-m+1)!} f^{(k)}(n_0+\theta_h h).
	\end{split}
\end{equation}
From \eqref{cn0fInd}, \eqref{cfkInd} and \eqref{eqah+1-ah}, we get
\begin{equation}
	\begin{split}
		a_{h+1} ~-~ a_h 
		&~<~  \sum_{i=1}^{k-m} \frac{(h+1)^i - h^i}{i!} \frac{1}{B_{m+i-1}} 
		~+~  
		2 \frac{(h+1)^{k-m+1}}{(k-m+1)!} \frac{1}{B_k}. 
		%\\
		%&~<~  \sum_{i=1}^{k-m} \frac{2^i h^{i-1}}{i!} \frac{1}{B_{m+i-1}} 
		%~+~  
		%2 \frac{(h+1)^{k-m+1}}{(k-m+1)!} \frac{1}{B_k}.
	\end{split}
\end{equation}
From \eqref{cBi6Ind}, we have
\begin{equation}\label{eqBm+i-1-Bm}
	B_{m+i-1} ~>~ B_m^{6^{i-1}} ~\geq~ B_m^{2(i-1)+1} 
	\phantom{m}\text{ and }~ 
	B_k ~>~ B_m^{6^{k-m}} ~\geq~ B_m^{2(k-m+1)+2}.
\end{equation}
Set
$$
R ~=~ R(m) ~=~ 2k! \Pi_H A_m.
$$
From \eqref{eqBm+i-1-Bm},
%and using the inequality $2^i \leq 2 \cdot i!$, 
for any $0 \leq h \leq 2R$, we get
\begin{equation}\label{eqlowah=1-ah}
	a_{h+1} ~-~ a_h 
	~<~ \frac{2}{B_m} \sum_{i=1}^{\infty}  \(\frac{2R}{B_m^2}\)^{i-1}
	~=~
	\frac{2}{B_m} \cdot \frac{1}{1-\frac{2R}{B_m^2}} 
	~<~ \frac{3}{B_m},
\end{equation}
since $\frac{2R}{B_m^2} < \frac{1}{4}$ by the choice of $R$ 
and \eqref{cB_kInd}. Also from \eqref{cn0fInd} and \eqref{eqah+1-ah}, we have
\begin{equation}\label{equppah+1-ah}
	a_{h+1} ~-~ a_h 
	~>~ \frac{1}{A_m} ~-~ 2 \frac{(4R)^{k-m+1}}{(k-m+1)!} \frac{1}{B_k} 
	~>~ \frac{1}{3A_m},
\end{equation}
for any $ 0 \leq h \leq 2R$. This is because, from \eqref{cB_kInd}, we get
$$
\frac{(4R)^{k-m+1}}{(k-m+1)!} \frac{1}{B_k} 
~<~ 2 \frac{(4k! \Pi_H A_m)^k}{(4k! \Pi_H A_{k-1})^{k+1}} 
~<~ \frac{1}{2k! \Pi_H A_{k-1}} < \frac{1}{3A_m}.
$$
Further, from \eqref{cn0fInd} and \eqref{cfkInd}, we have
\begin{equation*}
	\begin{split}
		a_R ~-~ a_0 
		&~=~ \sum_{i=1}^{k-m} \frac{R^i}{i!} \big\{ f^{(m+i-1)}(n_0) \big\} 
		~+~ \frac{R^{k-m+1}}{(k-m+1)!} f^{(k)}(n_0+\theta_R R) \\
		&~>~ \frac{R}{A_m} - \frac{R^{k-m+1}}{(k-m+1)!} \cdot \frac{1}{B_k} 
		~>~ 2k! \Pi_H -1,
	\end{split}
\end{equation*}
since from \eqref{cB_kInd},
$$
\frac{R^{k-m+1}}{(k-m+1)!} \cdot \frac{1}{B_k} 
~<~ \frac{(2k!\Pi_H A_m)^{k}}{(4 k!\Pi_H A_{k-1})^{k+1}} ~<~ 1.
$$
Hence there exists an integer $b \in (a_0, a_R)$ such that 
$$
b ~\equiv~ v_{m-1} ~(\mod k!\Pi_H).
$$
Let $r \in (0, R]$ be the least integer such that $a_r \geq b$. 
Then from \eqref{eqlowah=1-ah} and \eqref{equppah+1-ah}, we have
\begin{equation}\label{eqarintfr}
	\lf a_r \rf ~=~ b ~\equiv~ v_{m-1} ~(\mod k!\Pi_H) 
	\phantom{m}\text{and}\phantom{m}
	\{a_r\} ~<~ \frac{3}{B_m}.
\end{equation}
Set
$$
t ~=~ t(m) ~=~  
\( \bigg\lfloor \frac{4}{k! \Pi_H} \frac{A_m}{A_{m-1}} \bigg\rfloor +1 \) k! \Pi_H  
~-~ r_0,
$$
where $r_0 \in [0, k! \Pi_H)$ is an integer such that $r \equiv r_0 ~(\mod k! \Pi_H)$.
It is easy to note that $0 < t \leq R$.
Then from \eqref{eqarintfr} and
$$
a_{r+t} ~-~ a_r ~=~ \sum_{i = 1}^{t} \(a_{r+i} - a_{r+i-1} \) 
~<~ \frac{3t}{B_m} ~<~ \frac{15}{B_m} \frac{A_m}{A_{m-1}} 
~<~ \frac{15}{20} \frac{1}{B_{m-1}},
$$
we deduce that 
\begin{equation}\label{eqar+hintfracupp}
	\lf a_{r+t} \rf ~=~  \lf a_r \rf  ~\equiv~ v_{m-1} ~(\mod k!\Pi_H) 
	\phantom{m}\text{and}\phantom{m}
	\{a_{r+t}\} ~<~ \frac{1}{B_{m-1}}.
\end{equation}
Further, from \eqref{equppah+1-ah}, we have
$$
a_{r+t} ~-~ a_r ~>~ \frac{t}{3A_m} ~>~ \frac{1}{A_{m-1}}
$$
and hence we get
\begin{equation}\label{eqar+hfraclow}
	\{a_{r+t}\} ~>~ \frac{1}{A_{m-1}}.
\end{equation}
Set $s = r + t$ and $n_1=n_0 + s$. Observe that $s \equiv 0 ~(\mod k! \Pi_H)$ 
and $n_1 \in [n_0, n_0+2R]$. From \eqref{eqahdef}, \eqref{eqar+hintfracupp} and \eqref{eqar+hfraclow}, we deduce that
\begin{equation*}
	\lf f^{(m-1)}(n_1) \rf ~\equiv~ v_{m-1} ~(\mod k!\Pi_H) 
	\phantom{m}\text{and}\phantom{m}
	\{ f^{(m-1)}(n_1) \}  ~\in~ \(\frac{1}{A_{m-1}}, \frac{1}{B_{m-1}}\).
\end{equation*}
Now we will show that
\begin{equation*}
	\lf f^{(i)}(n_1) \rf ~\equiv~ v_i ~(\mod k!\Pi_H) 
	\phantom{m}\text{ and }\phantom{m}
	\{f^{(i)}(n_1)\} ~\in~ \(\frac{1}{2A_i}, \frac{2}{B_i}\)
\end{equation*}
for $m \leq i \leq k-1$.
Let $0 \leq j \leq k-m-1$. By Taylor's theorem, we have
\begin{equation}\label{eqTayf(m+j)}
	f^{(m+j)}(n_1) 
	~=~ f^{(m+j)}(n_0 +s) 
	~=~ \sum_{i=0}^{k-m-j-1} \frac{s^i}{i!} f^{(m+j+i)}(n_0) 
	~+~ \frac{s^{k-m-j}}{(k-m-j)!} f^{(k)}(n_0+ \theta s) 	
\end{equation}
for some $\theta \in [0,1]$. Note that from \eqref{cn0fInd} and \eqref{cB_kInd}, 
we get
\begin{equation}\label{eqfracf(m+j)low}
	\begin{split}
		\sum_{i = 0}^{k-m-j-1} \frac{s^i}{i!} \{ f^{(m+j+i)}(n_0) \}
		~+~ \frac{s^{k-m-j}}{(k-m-j)!} f^{(k)}(n_0+ \theta s)
		~>~ \frac{1}{A_{m+j}} -  \frac{(2R)^{k-m-j}}{(k-m-j)!} \frac{1}{B_k}
		~>~ \frac{1}{2A_{m+j}},
	\end{split}
\end{equation}
since
$$
\frac{(2R)^{k-m-j}}{(k-m-j)!} \frac{1}{B_k} 
~<~ 2 \frac{(2k!\Pi_H A_m)^{k}}{(4k!\Pi_H A_{k-1} )^{k+1}} 
~<~ \frac{1}{2A_{k-1}} ~<~ \frac{1}{2A_{m+j}}.
$$
Also, we have
\begin{equation*}
	\begin{split}
		\sum_{i = 0}^{k-m-j-1} \frac{s^i}{i!} \{ f^{(m+j+i)}(n_0) \}
		~+~ \frac{s^{k-m-j}}{(k-m-j)!} f^{(k)}(n_0 + \theta s) 
		~<~ \sum_{i = 0}^{k-m-j-1} \frac{s^i}{i!} \frac{1}{B_{m+j+i}} 
		~+~ \frac{s^{k-m-j}}{(k-m-j)!} \frac{1}{B_k}
	\end{split}
\end{equation*}
and from \eqref{cBi6Ind}, we get
\begin{equation*}
	\sum_{i = 0}^{k-m-j-1} \frac{s^i}{i!} \frac{1}{B_{m+j+i}} 
	~+~ \frac{s^{k-m-j}}{(k-m-j)!} \frac{1}{B_k}
	~<~ \frac{1}{B_{m+j}} \sum_{i=0}^{\infty} \(\frac{s}{B_{m+j}^2}\)^i
	~=~ \frac{1}{B_{m+j}} \cdot \frac{1}{1-\frac{s}{B_{m+j}^2}} 
	~<~ \frac{2}{B_{m+j}},
\end{equation*}
since
$$
\frac{s}{B_{m+j}^2} ~<~ \frac{2R}{B_{m}^2} ~<~ \frac{1}{2}.
$$
Thus we get
\begin{equation}\label{eqfracf(m+j)upp}
	\sum_{i = 0}^{k-m-j-1} \frac{s^i}{i!} \{ f^{(m+j+i)}(n_0) \}
	~+~ \frac{s^{k-m-j}}{(k-m-j)!} f^{(k)}(n_0+ \theta s) 
	~<~ \frac{2}{B_{m+j}}. 
\end{equation}
From \eqref{eqTayf(m+j)}, \eqref{eqfracf(m+j)low} and \eqref{eqfracf(m+j)upp}, 
we deduce that
$$
\lf f^{(m+j)}(n_1) \rf
~=~ \sum_{i=0}^{k-m-j-1} \frac{s^i}{i!} \lf f^{(m+j+i)}(n_0) \rf 
~\equiv~ v_{m+j} ~(\mod k!\Pi_H),
$$
since $s \equiv 0 ~(\mod k!\Pi_H)$ and
$$
\{ f^{(m+j)}(n_1) \} ~\in~ \(\frac{1}{2A_{m+j}}, \frac{2}{B_{m+j}}\)
$$
for $0 \leq j \leq k-m-1$. This completes the proof of \lemref{thmInd}.  \qed

\medspace

We now complete the proof of Theorem \ref{thmmain}. It is enough to show that if a function satisfies Conditions (\ref{zeroatinf}) and (\ref{plusatinf}), then it satisfies Conditions  \eqref{c0pr} - \eqref{c3pr} with $\ell = k-1$.

We consider an integer  $k \ge 2$ and a $k$-continuously differentiable real-valued function $f$ satisfying Relations (\ref{zeroatinf}) and (\ref{plusatinf}). We let $H$ be an integer which is larger than or equal to $k$; we recall that $\Pi_H$ denote the product of the primes up to $H$ and notice the inequality $\Pi_H \ge H$ for $H \ge 2$.

Set
\begin{equation}\label{chBi}
\begin{split}
D_0 ~=~ 4 \cdot 2^k, \phantom{mm} 
D_1 ~=~ 2^{8k} k! \Pi_H, \phantom{mm} 
D_{i+1} ~=~ 2 D_i^6 \phantom{m}\text{for}\phantom{m} 1 \leq i \leq k-2, 
\end{split}
\end{equation}
\begin{equation}\label{chAiBk}
C_i ~=~ 2^{5(k-i)} D_i \phantom{m} \text{for}\phantom{m} 0 \leq i \leq k-1 \phantom{mm}\text{and}\phantom{mm}
D_k ~=~ \(4 k! \Pi_H 2^k C_{k-1}\)^{6(k+1)}.
\end{equation}
From \eqref{zeroatinf}, there exists an integer $x_0 > 1$ such that
\begin{equation*}
|f^{(k)}(x)| ~<~ \frac{1}{D_k} \phantom{m}\text{for}\phantom{m} x \geq x_0.
\end{equation*}

By (\ref{plusatinf}) and the continuity of $f^{(k-1)}$, any sufficiently large real number is a value of $f^{(k-1)}$ and thus one can find a positive integer $s$ and a real $x_1 > x_0$ such that
$$
f^{(k-1)} (x_1) = (k!)^s\Pi_H + \frac{1}{2D_{k-1}}.
$$
We let $n_{k-1} = \lf x_1 \rf +1$. By the mean value theorem, we have
\begin{equation}\label{eqfk-1fracint}
	\frac{1}{C_{k-1}} \le \frac{1}{2D_{k-1}}-\frac{1}{D_k} < f^{(k-1)}(n_{k-1})-(k!)^s\Pi_H <  \frac{1}{2D_{k-1}}+\frac{1}{D_k} \le \frac{1}{D_{k-1}}.
\end{equation}
Thus we have
\begin{equation}\label{eqindr=1}
\lf f^{(k-1)}(n_{k-1}) \rf ~=~ (k!)^s\Pi_H 
\phantom{mm}\text{and}\phantom{mm}
\{  f^{(k-1)}(n_{k-1}) \} ~\in~ \(  \frac{1}{C_{k-1}},~  \frac{1}{D_{k-1}}  \).
\end{equation}
Set
$$
v_0 ~=~ 1  
\phantom{mm}\text{and}\phantom{mm}
v_i ~=~ 0 \phantom{m}\text{for}\phantom{m} 1 \leq i \leq k-1.
$$
We inductively show that 
for each integer $1 \leq r \leq k$, there exists an integer 
$n_{k-r} \in [n_{k-1}, n_{k-1}+ 4k!\Pi_H 2^r C_{k-1}]$ such that
\begin{equation}\label{eqIndHypfi}
	\lf f^{(i)}(n_{k-r}) \rf ~\equiv~ v_i ~(\mod k! \Pi_H) 
	\phantom{m}\text{and}\phantom{m}
	\{f^{(i)}(n_{k-r})\} ~\in~ \(\frac{1}{2^{r-1} C_i}, \frac{2^{r-1}}{D_i}\) 
	\phantom{m}
	\text{ for } k-r \leq i \leq k-1 .
\end{equation}
From \eqref{eqindr=1}, the statement \eqref{eqIndHypfi} is true when $r=1$. Now suppose that the statement \eqref{eqIndHypfi} is true for $r = t\in \{1, 2, \cdots, k-1\}$. As a consequence of \lemref{thmInd}, we will show that the statement  \eqref{eqIndHypfi} is true for $r = t+1$. Set
\begin{equation}
	A_{i} ~=~ 2^{t-1} C_i 
	\phantom{m} \phantom{m}\text{for } k-t-1 \leq i \leq k-1
	\phantom{m}\text{and}\phantom{m}
	B_{i} ~=~ \frac{D_i}{2^{t-1}}
	\phantom{m}\text{for } k-t-1 \leq i \leq k.
\end{equation}
By induction assumption, there exists an integer 
$n_{k-t} \in [n_{k-1},  n_{k-1}+ 4k!\Pi_H 2^t C_{k-1}]$ such that
\begin{equation}
	\lf f^{(i)}(n_{k-t}) \rf ~\equiv~ v_i ~(\mod k! \Pi_H) 
	\phantom{m}\text{and}\phantom{m}
	\{f^{(i)}(n_{k-t})\} ~\in~ \(\frac{1}{A_i},~ \frac{1}{B_i}\) \phantom{m}
	\text{ for } k-t \leq i \leq k-1 .
\end{equation}
It is easy to check that the conditions \eqref{cn0fInd} - \eqref{cAm-1Bm-1Ind} in \lemref{thmInd} are satisfied with $m = k-t$ and $n_{k-t}$ in place of $n_0$.
Hence by applying \lemref{thmInd}, we deduce that there exists an integer 
$n_{k-t-1} \in [n_{k-t}, n_{k-t} + 4 k!\Pi_H A_{k-t}]$ such that
\begin{equation*}
	\lf f^{(i)}(n_{k-t-1}) \rf ~\equiv~ v_i ~(\mod k! \Pi_H) 
	\phantom{m}\text{and}\phantom{m}
	\{f^{(i)}(n_{k-t-1})\} ~\in~ \(\frac{1}{2A_i},~ \frac{2}{B_i}\) 
	\phantom{m}
	\text{ for } k-t-1 \leq i \leq k-1 .
\end{equation*}
Thus we have
$$
\{f^{(i)}(n_{k-t-1})\} ~\in~ \(\frac{1}{2^t C_i},~ \frac{2^t}{D_i}\) \phantom{m}
\text{ for } k-t-1 \leq i \leq k-1.
$$
It is easy to see that $n_{k-t-1 } \in [n_{k-1}, n_{k-1} +  4 k!\Pi_H 2^{t+1} C_{k-1}]$.
Thus the statement \eqref{eqIndHypfi} is true for $r = t+1$. Hence 
we conclude by induction that the  statement \eqref{eqIndHypfi} is true for each $1 \leq r \leq k$. Thus there exists an integer 
$n_{0} \in [n_{k-1}, n_{k-1}+ 4k!\Pi_H 2^{k} C_{k-1}]$ 
such that
\begin{equation}\label{eqintcong}
	\lf f^{(i)}(n_{0}) \rf ~\equiv~ v_i ~(\mod k! \Pi_H) 
	\phantom{m}\text{and}\phantom{m}
	\{f^{(i)}(n_{0})\} ~\in~ \(\frac{1}{2^{k-1} C_i}, \frac{2^{k-1}}{D_i}\) 
	\phantom{m}
	\text{ for } 0 \leq i \leq k-1 .
\end{equation}
Then we have
\begin{equation}\label{eqfinosfrc}
\{f(n_{0})\} ~\in~ \(\frac{1}{4^{4k}},~ \frac{1}{4}\)
\phantom{mm}\text{and}\phantom{mm}
\{f^{(i)}(n_{0})\} ~<~ \frac{1}{4H^i}, ~~1 \leq i \leq k-1,
\end{equation}
since 
$$
D_i~>~ D_1^{6^{i -1}} ~=~ \(2^{8k} k! \Pi_H\)^{6^{i -1}} ~>~ 4H^i \phantom{m}
\text{ for }\phantom{m} 1 \leq i \leq k-1 .
$$
Also we have
\begin{equation}\label{eqfkupp}
|f^{(k)}(x)| ~<~ \frac{1}{D_k} ~<~ \frac{1}{4^{6k} H^k} \phantom{m}
\text{ for }\phantom{m} x \geq n_{0}.
\end{equation}
By the mean value theorem, we get
\begin{equation}\label{eqmvtdiffnonk-1}
|f^{(k-1)}(n_{0}) - f^{(k-1)}(n_{k-1})| 
~<~ \frac{n_{0} - n_{k-1}}{D_k} 
~\leq~ \frac{ 4k!\Pi_H 2^k C_{k-1} }{ \(4 k! \Pi_H 2^k C_{k-1}\)^{6(k+1)} } 
~<~ \frac{1}{C_{k-1}}.
\end{equation}
From \eqref{eqfk-1fracint}, \eqref{eqmvtdiffnonk-1} and by noting
$$
\frac{1}{D_{k-1}} + \frac{1}{C_{k-1}} < \frac{2}{D_{k-1}} < 1,
$$
we deduce that
\begin{equation}\label{eqfk-1int}
\lf f^{(k-1)}(n_{0}) \rf ~=~ (k!)^s \Pi_H.
\end{equation}
From \eqref{eqintcong} and \eqref{eqfk-1int}, we get
\begin{equation}\label{eqgcdf0k-1nos}
\gcd\(\lf f(n_{0}) \rf ,~ \lf f^{(k-1)}(n_{0}) \rf\) ~=~ 1.
\end{equation}
Thus it follows from \eqref{eqintcong}, \eqref{eqfinosfrc}, \eqref{eqfkupp} and \eqref{eqgcdf0k-1nos} that there exist infinitely many positive integers $n$ satisfying \eqref{c0pr} - \eqref{c3pr}. 
Now \thmref{thmmain} follows from  \propref{reduction}. \qed

\medspace

\section{Proof of \thmref{thmB}}
Let the notations be as before. Let $H_1 \geq k$ be a positive integer. From Section \ref{Sthmmain}, there exists a positive integer $n_1$ satisfying \eqref{c0pr} - \eqref{c3pr} with $H = H_1$. Hence from \propref{reduction}, we have
$$
\gcd\( \lf f(n_1 + i_1)\rf,  \lf f(n_1 + i_2)\rf, \cdots,  \lf f(n_1 + i_k)\rf\) ~=~ 1 
$$
for any integers $1 \leq i_1 < i_2 < \cdots < i_k \leq H_1$. Let $H_2$ be a positive integer such that 
$$
H_2 ~>~ H_1 + \sum_{j = 1}^{H_1} |f(n_1 + j)|.
$$
Again from Section \ref{Sthmmain}, there exists a positive integer $n_2 > n_1 + H_1$ satisfying \eqref{c0pr} - \eqref{c3pr} with $H = H_2$. From Section \ref{secreduction}, we have $p \nmid \lf f(n_2 + h_2)\rf$ for any prime $p \leq H_2$ and any positive integer $h_2 \leq H_2$.
Hence $\gcd\(\lf f(n_1 + h_1)\rf, \lf f(n_2 + h_2)\rf\) = 1$ for any positive integers $h_1 \leq H_1$ and $h_2 \leq H_2$. Inductively, we can show that for a positive integer $H_r$ satisfying
$$
H_r ~>~ H_{r-1} + \sum_{j = 1}^{H_{r-1}} |f(n_{r-1} + j)|,
$$
there exists a positive integer $n_r > n_{r-1} + H_{r-1}$ satisfying \eqref{c0pr} - \eqref{c3pr} with $H = H_r$. Set
$$
\cA ~=~ \{n_i + h_i ~:~ h_i \in [1, H_i] \cap \N,~ i \in \N\}.
$$
It is easy to see that the set $\cA$ has upper Banach density equal to $1$ and 
$$
\gcd\( \lf f(m_1)\rf,  \lf f(m_2)\rf, \cdots,  \lf f(m_k)\rf\) ~=~ 1 
$$
for any distinct integers $m_1, m_2, \cdots, m_k \in \cA$. This completes the proof of \thmref{thmB}. \qed

\medspace

\section{Acknowledgments}
The authors would like to acknowledge SPARC project 445. When preparing this paper, the first-named author benefitted from the support of the joint FWF-ANR project ArithRand: FWF: I 4945-N and ANR-20-CE91-0006. 
The second author would like to thank Queen's University, Canada, for providing an excellent atmosphere for work. The second author would also like to acknowledge the Institute of Mathematical Sciences (IMSc), India, and the Institut de Math\'ematiques de Bordeaux (IMB), France, where this work was initiated.

\end{document}